\documentclass[12pt,a4paper,reqno]{amsart}

\usepackage{amsfonts}
\usepackage{amsmath}
\usepackage{amssymb}
\usepackage{amsthm}
\usepackage{amsxtra}
\usepackage{bbm}
\usepackage{braket}
\usepackage{comment}
\usepackage{extarrows}
\usepackage{graphicx}
\usepackage{latexsym}
\usepackage{mathrsfs}
\usepackage{yhmath}
\usepackage{nicematrix}
\usepackage{mathtools}

\usepackage{float}
\usepackage{booktabs}
\usepackage{verbatim}
\usepackage{xcolor}

\usepackage[pdfborder={0 0 0}]{hyperref} % hyperref must be loaded after natbib
\hypersetup{breaklinks, raiselinks}

\usepackage{bm}
\usepackage{tikz}
\usepackage[initials,lite]{amsrefs} % amsrefs must be loaded after hyperref
\AtBeginDocument{%
    \def\MR#1{}
}
\usepackage{cleveref}

\usepackage{a4wide}
\usepackage{fullpage}

\theoremstyle{plain}
\newtheorem{Theorem}{Theorem}[section]
\newtheorem{Lemma}[Theorem]{Lemma}
\newtheorem{Corollary}[Theorem]{Corollary}
\newtheorem{Proposition}[Theorem]{Proposition}

\theoremstyle{definition}

\newtheorem{Assumptions and Discussion}[Theorem]{Assumptions and Discussion}

\newtheorem{Example}[Theorem]{Example}
\newtheorem{Definition}[Theorem]{Definition}

\newtheorem{Question}[Theorem]{Question}

\newtheorem{Remark}[Theorem]{Remark}

\newtheorem{Notation}[Theorem]{Notation}

\theoremstyle{remark}

\newtheorem{Setting}[Theorem]{Setting}

\newtheorem*{acknowledgment*}{Acknowledgment}

\usepackage[shortlabels]{enumitem}
\SetEnumerateShortLabel{a}{\textup{(\alph*)}}
\SetEnumerateShortLabel{A}{\textup{(\Alph*)}}
\SetEnumerateShortLabel{1}{\textup{(\arabic*)}}
\SetEnumerateShortLabel{i}{\textup{(\roman*)}}
\SetEnumerateShortLabel{I}{\textup{(\Roman*)}}
\def\lex{\operatorname{lex}}

\def\deg{\operatorname{deg}}
\def\dim{\operatorname{dim}}

\def\floor#1{\left\lfloor #1 \right\rfloor}

\def\Tor{\operatorname{Tor}}

\def\id{\operatorname{id}}

\def\KK{{\mathbb K}}
\def\lex{{\operatorname{lex}}}

\def\Pol{\operatorname{Pol}} % polarization
\def\PP{{\mathbb P}}

\def\suppdeg{\operatorname{suppdeg}}

\def\type{\operatorname{type}}

\def\ZZ{{\mathbb Z}}

\newcommand\bda{{\bm a}}

\newcommand\bdb{{\bm b}}

\newcommand\bdF{{\bm F}}

\newcommand\bdG{{\bm G}}

\newcommand\bdlambda{{\bm \lambda}}

\newcommand\bdmu{{\bm \mu}}

\newcommand\bdZ{{\bm Z}}
\newcommand\bdz{{\bm z}}

\newcommand\calF{\mathcal{F}}
\newcommand\calG{\mathcal{G}}

\newcommand\calI{\mathcal{I}}
\newcommand\calJ{\mathcal{J}}

\newcommand\calP{\mathcal{P}}

\newcommand\frakm{\mathfrak{m}}

\newcommand\frakp{\mathfrak{p}}

\def\frakS{\mathfrak{S}}

\newcommand{\Ass}{\operatorname{Ass}}

\newcommand{\pd}{\operatorname{pd}}

\begin{document}

\title[Generalized Star Configurations]{Symbolic powers of generalized star configurations of hypersurfaces}

\author[Kuei-Nuan Lin, Yi-Huang Shen]{Kuei-Nuan Lin and Yi-Huang Shen}

%\thanks{\today}

\thanks{2020 {\em Mathematics Subject Classification}.
    %    Primary 13F55, 13P10, 14M25, 16S37; Secondary 14M05, 14N10, 05E40, 05E45
    13A15, %	Ideals; multiplicative ideal theory 
    13A50, %  	Actions of groups on commutative rings; invariant theory
    13D02, % 	Syzygies, resolutions, complexes
    14N20, % 	Configurations and arrangements of linear subspaces
    52C35. %    Arrangements of points, flats, hyperplanes
}


\address{Department of Mathematics, The Penn State University, McKeesport, PA, 15132, USA}
\email{kul20@psu.edu}

\address{CAS Wu Wen-Tsun Key Laboratory of Mathematics, School of Mathematical Sciences, University of Science and Technology of China, Hefei, Anhui, 230026, P.R.~China}
\email{yhshen@ustc.edu.cn}

\begin{abstract}
    We introduce the class of sparse symmetric shifted monomial ideals. These
    ideals have linear quotients and their Betti numbers are computed. Using
    this, we prove that the symbolic powers of the generalized star
    configuration ideal are sequentially Cohen--Macaulay under some mild
    genericness assumption. With respect to these symbolic powers, we also
    consider the Harbourne--Huneke containment problem and establish the
    Demailly-like bound.
\end{abstract}

\maketitle
\section{Introduction}

In the simplest case, a \emph{star configuration} of codimension $c$ in $\PP^n$
is a certain union of linear subspaces $V_1,\dots,V_p$ each of codimension $c$.
The study of various aspects of generalizations of star configurations has
attracted considerable attention recently; see for example, \cite{MR2492444},
\cite{MR3197149}, \cite{MR3003727}, \cite{MR3683102}, \cite{arXiv:1912.04448},
\cite{MR4177276} and \cite{arXiv:1906.08346}, to name just a few.  As a sequel
to our previous work \cite{arXiv:1912.04448} along this line, we continue to
investigate the symbolic powers of the generalized star configuration of
hypersurfaces. 

Let $R=\KK[x_0,\dots,x_n]$ be a standard graded polynomial ring over an
infinite $\KK$ and denote its graded maximal ideal by $\frakm$. Suppose that
$s\ge n+1$ and let $\calF=\{f_1,\dots,f_s\}$ be a set of forms in $R$
satisfying some genericness condition.  For a fixed choice of positive integers
$a$ and $b$ such that $a\le bs$, the \emph{generalized star configuration
    ideal} that we shall study here is the \emph{uniform $a$-fold product
    ideal} $I_a(\calF^b)\coloneqq I_a(f_1^b\cdots f_s^b)$.

In general, the \emph{$a$-fold product ideal} $I_a(f_1^{b_1}\cdots f_s^{b_s})$
is generated by the $a$-fold products of the forms $f_1,\dots,f_s$ with
multiplicities $b_1,\dots,b_s$ respectively:
\begin{equation*}
    \Braket{f_1^{n_1}\cdots f_s^{n_s} : 0\le n_i\le b_i\text{ for each $i$ such that $\sum\nolimits_{i}n_i=a$}}.
    %\label{eqndef:a-fold-product-ideal}
\end{equation*}
It has also been studied extensively recently; see for example \cite{zbMATH06759435} and \cite{MR3864202}. This notion was originally introduced as a nice tool for determining the minimum distance of linear codes in the coding theory. And it also emerges naturally when dealing with the higher order Orlik--Terao algebra of hyperplane arrangements.
%\alert{backgroud like coding theory}
%And even for general $b$, the ideal $I_a(\calF^b)$ still exhibits many similar properties by our previous work in \cite{arXiv:1912.04448}.

The $a$-fold product ideal $I_a(f_1^{b_1}\cdots f_s^{b_s})$ is \emph{uniform} when the multiplicities $b_1,\dots,b_s$ are all identical, say, to some positive integer $b$. Whence, we will write it as $I_a(\calF^b)$ for simplicity. It is treated as the \emph{generalized star configuration ideal}, since when $b=1$, we will get back the star configuration ideal of hypersurfaces 
\[
    I_{c,\calF} \coloneqq \bigcap_{1\le i_1<\cdots<i_c\le s}\braket{f_{i_1},\dots,f_{i_c}} 
\]
for $c=s-a+1$, which were studied for instance in \cite{MR3683102} and \cite{MR4177276}.  
In their work, the set $\calF$ is often assumed to be \emph{$c$-generic}, namely any subset of size at most $c+1$ will form a regular sequence.

Our current work will focus on the study of the symbolic powers of the generalized star configuration ideal $I_a(\calF^b)$.
Recall that give an ideal $I$ in a ring $R$, the \emph{$m$-th symbolic power} of $I$ is defined by
\[
    I^{(m)}\coloneqq \bigcap_{\frakp\in \Ass(R/I)} I^mR_\frakp \cap R.
\]
For instance, when $I$ is the defining ideal of a reduced affine scheme over an algebraically closed field of characteristic zero, Zariski and Nagata showed that $I^{(m)}$ is generated by the polynomials whose partial derivatives of orders up to $m-1$ vanish on this scheme.
%the contraction of $I^mR_W$ to $R$ where $W$ is the complement of the union of the associate primes of $I$ and $R_W$ is the localization of $R$ at the multiplicative system $W$.
%Recall that the $m$-th symbolic power of this ideal can be identified as the one of all hypersurfaces in $\mathbb{P}^n_{\mathbb{K}}$ vanishing at $X$ with order at least $m$. 
And in the previous work \cite{arXiv:1912.04448}, we have already studied the resurgence and symbolic defect of $I_a(\calF^b)$, as well as the Betti table and the primary decomposition of its symbolic powers. %the $a$-fold product of the generalized star configuration \cite{arXiv:1912.04448}.  

In \cite{MR3003727} and \cite{MR3683102},  Geramita, Harbourne, Migliore, and Nagel showed that the symbolic powers of the ordinary star configuration ideal define arithmetically Cohen–Macaulay subschemes. It is then natural to ask if the same property still holds in the generalized case. 
Notice that the ideal $I_a(\calF^b)$ is no longer unmixed for $b\ge 2$. Thus, the correct question can be put as: whether the symbolic powers of this ideal are all sequentially Cohen--Macaulay?
We provide a positive answer under some mild genericness assumption.

In this work, we also study the Harbourne--Huneke containment problem, which was originally raised in \cite{MR3115195}. After \cite{arXiv:2009.05022}, a generalized version can be stated as: given a homogeneous ideal $I$ of big height $h$ in a standard graded ring with the maximal homogeneous ideal $\mathfrak{m}$, does the inequality 
\[
    I^{(\ell(h+m-1))}\subseteq \frakm^{\ell(h-1)}(I^{(m)})^\ell
\]
hold for all $m,\ell\ge 1$?
The comparison between ordinary powers of an ideal and its symbolic ideals is wildly open in general; see, for example, \cite{MR2629595} and \cite{MR1881923}. We provide a positive answer for Harbourne--Huneke containment problem in the case $I=I_a(\calF^b)$ 
%(\Cref{thmHH})
under some mild genericness assumption. As a by-product, we establish the Demailly-like bound 
%(\Cref{Delike})
in this case as well.

%Here is the plan of this paper.
%In order to show the sequentially Cohen--Macaulayness of the uniform $a$-fold product ideals, we work on the monomial case first, i.e., when the hypersurfaces are defined by the variables of the ring. We will introduce the notion of \emph{sparse symmetric shifted ideals}. This class of squarefree monomial ideal is a generalization of the symmetric shifted ideals introduced in \cite{MR4108337}. Sparse symmetric shifted  ideals give similar properties as in symmetric shifted case. We show it has linear quotients (\Cref{SSSisLQ}), and give its Betti numbers (\Cref{SSSBetti}).

As usual, for a homogeneous ideal $J$, we will write $\alpha(J)$ for the least degree of nonzero forms in this ideal. Then, our main results can be summarized in the following.
\begin{Theorem}  
    Let $a$ and $b$ be positive integers such that $(b-1)s+1<a\le bs$.  We assume further that  $n\ge  bs-a+1$ and $\calF=\{f_1,\dots,f_s\}$ is a set of $(bs-a+1)$-generic $d$-forms in $R=\KK[x_0,\dots,x_n]$. Then the following properties hold for the generalized star configuration ideal $I=I_a(\calF^b)$ whose big height is known to be $h=bs-a+1$.
    \begin{enumerate}[a]
        \item The symbolic powers of $I$ are all sequentially Cohen--Macaulay.
        \item The ideal $I$ satisfies
            \begin{equation*}
                I^{(\ell(h+m-1)-h+k)}\subseteq \frakm^{d((\ell-1)(h-1)+k-1)(bh-(bs-a))}(I^{(m)})^\ell
                %\label{eqn-5}
            \end{equation*}
            for all positive integers $k,\ell$, and $m$.
        \item The inequality 
            \[
                \frac{\alpha(I^{(\ell)})}{\ell}\ge \frac{\alpha(I^{(m)})+h-1}{m+h-1}
            \]
            holds for all $\ell, m \ge 1$.     
    \end{enumerate}
\end{Theorem}

% \begin{Theorem}  
%     Let $a,b$ be positive integers such that $(b-1)s+1<a\le bs$ and set $\delta\coloneqq a-(b-1)s-1$.  We assume further that  $\delta\ge s-n$ and $\calF=\{f_1,\dots,f_s\}$ is a set of $(s-\delta)$-generic $d$-forms in $R=\KK[x_0,\dots,x_n]$. Then the following properties hold for the generalized star configuration ideal $I=I_a(\calF^b)$ whose big height will be denoted by $h$.
%     \begin{enumerate}[a]
%         \item The symbolic powers of $I$ are all sequentially Cohen--Macaulay.
%         \item The ideal $I$ satisfies
%             \begin{equation*}
%                 I^{(\ell(h+m-1)-h+k)}\subseteq \frakm^{d((\ell-1)(h-1)+k-1)(\mu_a^0+b(h-c_0))}(I^{(m)})^\ell
%                 %\label{eqn-5}
%             \end{equation*}
%             for all positive integers $k,\ell$ and $m$.
%         \item The inequality 
%             \[
%                 \frac{\alpha(I^{(\ell)})}{\ell}\ge \frac{\alpha(I^{(m)})+h-1}{m+h-1}
%             \]
%             holds for all $\ell, m \ge 1$.     
%     \end{enumerate}
% \end{Theorem}

\section{Sparse symmetric shifted ideals}
%\label{section3}

To show the sequentially Cohen--Macaulayness of the symbolic powers of the generalized star configuration ideal $I_a(\calF^b)$, we have to introduce first the sparse symmetric shifted ideals. These squarefree monomial ideals can be viewed as variations of the symmetric shifted ideals introduced in \cite{MR4108337}. Likewise, the new ideals have linear quotients (\Cref{SSSisLQ}) and the Betti tables can be obtained combinatorially with ease (\Cref{SSSBetti}).

\subsection{Definition and the linear quotient property}
The ground ring that we shall consider is the polynomial ring $Z=\KK[z_{i,j}:1\le i\le s, 1\le j\le b]$ over an arbitrary field $\KK$. 
For each tuple $\bdmu=(\mu_1,\dots,\mu_s)\in \{0,1,2,\dots,b\}^s$,
a monomial of the form $\bdz^{\bdmu}\coloneqq \prod_{i=1}^s z_{i,\mu_i}$ will be called a \emph{sparse squarefree monomial}. Here,  when $\mu_i=0$, $z_{i,\mu_i}$ is understood to be the identity element $1\in Z$.

Recall that a tuple $\bdlambda=(\lambda_1,\dots,\lambda_s)\in \ZZ_{\ge 0}^s$ is called a \emph{partition} of length $s$, if $\lambda_1\le \cdots \le \lambda_s$. For the sparse squarefree monomial $f=\bdz^{\bdmu}$ above, we may suppose that the multi-set
\[
    \{\mu_1,\dots,\mu_s\}=
    \{\lambda_1\le \cdots \le \lambda_s\}.
\]
Now, we write
\[
    \bdlambda(f)\coloneqq (\lambda_1,\dots,\lambda_s)\in \ZZ_{\ge 0}^s,
\]
and call it the \emph{partition} associated to $f$. 
If $\bdlambda(f)=(0,\dots,0,e_1,\dots,e_d)$ where $e_1\ge 1$, then
the \emph{support degree} of both $f$ and $\bdlambda=\bdlambda(f)$ are 
\[
    \suppdeg(f)=\suppdeg(\bdlambda)\coloneqq d,
\]
and the \emph{weight} of both $f$ and $\bdlambda$ are 
\[
    w(f)=|\bdlambda|\coloneqq e_1+\cdots+e_d=\lambda_1+\cdots+\lambda_s=\mu_1+\cdots+\mu_s.
\]
Related, a sparse squarefree monomial in $Z$ whose support degree is $s$, is called a \emph{rainbow monomial} in \cite[Subsection 1.3]{arXiv.1912.03898}.

Let $\bdmu'=(\mu_1',\dots,\mu_s')\in \{0,1,2,\dots,b\}^s$ be another tuple. We will write $\bdmu'\le_0 \bdmu$ if
\begin{itemize} 
    \item $\mu_i'\le \mu_i$ for each $i$, and
    \item $\mu_i'=0$ if and only if $\mu_i=0$.
\end{itemize}

Notice that the permutation group $\frakS_s$ has an action on $Z$, such that for each permutation $\sigma\in \frakS_s$ and a variable $z_{i,j}$ in $Z$, we have $\sigma(z_{i,j})=z_{\sigma(i),j}$.  

\begin{Definition}
    A squarefree monomial proper ideal $J$ of $Z$ will be called \emph{sparse symmetric shifted} if it satisfies the following requirements.
    \begin{enumerate}[1]
        \item[(S-1)] The ideal $J$ is generated by some sparse squarefree monomials.
        \item[(S-2)] For each minimal monomial generator $f=\bdz^{\bdmu}\in G(J)$ with $\bdmu=(\mu_1,\dots,\mu_s)$, if $f'=\bdz^{\bdmu'}$ with $\bdmu'=(\mu_1',\dots,\mu_s')$ and $\bdmu'\le_0 \bdmu$, then $f'\in J$.
        \item[(S-3)] The ideal $J$ is invariant under the action of $\frakS_s$.
    \end{enumerate}
\end{Definition}

\begin{Remark}
    \label{rmk:sss}
    We may imagine the variables in $Z$ displayed in the canonical matrix $\bdZ=(z_{i,j})_{s\times b}$. Then, roughly speaking, the sparse symmetry shifted ideal above has sparseness, vertical (upward) shiftiness, and horizontal symmetry. And these properties correspond to the (S-1), (S-2), and (S-3) respectively. Meanwhile, the symmetric shifted ideal, introduced in \cite{MR4108337}, has both horizontal shiftiness and horizontal symmetry, since it lacks the sparse structure here.
\end{Remark}

\begin{Lemma}
    \label{lem-1}
    Let $J$ be a sparse symmetric shifted ideal in $Z$. Suppose that there is some sparse squarefree monomial $f=\prod_i z_{\sigma(i),\mu_i}\in J$ for some $\sigma\in \frakS_s$. If $f'=\prod_i z_{\sigma'(i),\mu_i'}$ for some $\sigma'\in \frakS_s$ with $(\mu_1',\dots,\mu_s')\le_0 (\mu_1,\dots,\mu_s)$, then $f'\in J$ as well.
\end{Lemma}

\begin{proof}
    By the symmetry condition (S-3), it suffices to consider the case when $\sigma=\sigma'=\id\in \frakS_s$.  Since $f\in J$, there exists some minimal monomial generator $g=\prod_i z_{i,\tau_i}\in G(J)$ such that $\tau_i=\mu_i$ or $0$ for each $i$. Define $\tau_i'=\min(\tau_i,\mu_i')$ for each $i$. Since $(\tau_1',\dots,\tau_s')\le_0 (\tau_1,\dots,\tau_s)$, the monomial $g'\coloneqq \prod_i z_{i,\tau_i'}\in J$ by the shiftiness condition (S-2). As $g'$ divides $f'$, this implies that $f'\in J$ as well. 
\end{proof}

\begin{Corollary}
    \label{lem-2}
    Let $J$ be a sparse symmetric shifted ideal in $Z$. Suppose that $f=\prod_{i}z_{\sigma(i),\lambda_i}\in J$ where $\bdlambda(f)=(\lambda_1,\dots,\lambda_s)$ and $\sigma\in \frakS_s$. If $\lambda_{i_1}=0<k\le\lambda_{i_2}$ for some $i_1$ and $i_2$, then $f'\coloneqq f z_{\sigma(i_1),k}/z_{\sigma(i_2),\lambda_{i_2}}\in J$.
\end{Corollary}

%\begin{proof}
%    Consider $f''\coloneqq fz_{\sigma(i),\lambda_j}/z_{\sigma(j),\lambda_j}$. Since $\bdlambda(f'')=\bdlambda(f)$, it follows from (S-3) that the sparse squarefree monomial $f''\in J$. Now we apply \Cref{lem-1} to the pair $(f'', f')$.
%\end{proof}

\begin{Remark}
    \label{rmk:colon}
    Suppose that the ideal $J\subset Z$ is generated by some sparse squarefree monomials. If $f$ is a monomial in $Z$, then the colon ideal $J:_Z f$ is also generated by some sparse squarefree monomials that are coprime to $f$. To see this, it suffices to notice that when $g$ is a sparse squarefree monomial, then $\braket{g}:f$ is generated by the sparse squarefree monomial $g/\gcd(f,g)$.
\end{Remark}

Here is the first main result of this section.

\begin{Theorem}
    \label{SSSisLQ}
    If $J$ is a sparse symmetric shifted ideal in $Z$, then it has linear quotients.
\end{Theorem}

%The next two subsections are devoted to proving this result.
\begin{proof}
    Inspired by the treatment in \cite{MR4108337}, we first consider the lexicographic ordering $>_{\lex}$ on $Z$ where the variables $z_{i,j}>z_{i',j'}$ if and only if
    \begin{itemize}
        \item $j<j'$, or
        \item $j=j'$ and $i<i'$.
    \end{itemize}
    Thus, among those squarefree monomials of degree $s$, $z_{1,1}z_{2,1}\cdots z_{s,1}$ is the initial (biggest) one with respect to $>_{\lex}$.

    Next, we take a total order $\prec$ on $\ZZ_{\ge 0}^s$ defined by
    \[
        \bda=(a_1,\dots,a_s) \prec \bdb=(b_1,\dots,b_s)
    \]
    if and only if
    \begin{itemize}
        \item $\suppdeg(\bda)<\suppdeg(\bdb)$, or
        \item $\suppdeg(\bda)=\suppdeg(\bdb)$ and $|\bda|<|\bdb|$, or
        \item $\suppdeg(\bda)=\suppdeg(\bdb)$ and $|\bda|=|\bdb|$, while the leftmost non-zero entry of $(a_1 - b_1, \dots, a_s - b_s)$ is negative. 
    \end{itemize} 

    Now, we are ready to introduce a total order on the minimal generating set $G(J)$, which, by abuse of notation, will still be denoted by $\prec$. For $f_1$ and $f_2$ in $G(J)$, we define $f_1\prec f_2$ if and only if
    \begin{itemize}
        \item $\bdlambda(f_1)\prec\bdlambda(f_2)$, or
        \item $\bdlambda(f_1)=\bdlambda(f_2)$ and $f_1>_{\lex}f_2$. 
    \end{itemize}
    Then, we sort the squarefree monomials in $G(J)$ with respect to $\prec$. To establish the linear quotient property, we take arbitrary minimal monomial generator $f\in G(J)$ and let 
    \[
        J_f\coloneqq \braket{g\in G(J):g\prec f}.
    \]
    The remaining task is to show that the colon ideal $J_f:_{Z} f$ is linear for each such $f$. 
    We can accomplish this task with the following two steps. 
    \begin{enumerate}[i]
        \item First of all,  we describe explicitly the set $\calG_f$ of variables that belong to this colon ideal.
        \item Next, we prove that the set $\calG_f$ actually generates the colon ideal $J_f: f$ in general. Thus, $J_f:f$ is linear, as expected.
    \end{enumerate}
    They are established in the coming \Cref{clm:G} and \Cref{lem:containment}
    %We will prove the above two steps in the coming two subsections 
    respectively. And consequently, the ideal $J$ has linear quotients.
\end{proof}

%\subsection{Linear generators}
%\label{ssLG}
In the remaining of this subsection, we fix a monomial $f$ in $G(J)$. The first task is to describe the set $\calG_f$ which collects the variables in the colon ideal $J_f: f$.  
Suppose that the partition $\bdlambda(f)=(\lambda_1,\dots,\lambda_s)$. Now, there is some permutation $\sigma\in \frakS_s$ such that $f=\prod_{i=1}^s z_{\sigma(i),\lambda_i}$. Recall that for $i$ with $\lambda_i=0$, $z_{\sigma(i),\lambda_i}$ is understood to be the identity element $1\in Z$. 
%Now, we want to decide whether $z_{\sigma(i),j}$ belongs to $G(J_f:f)$, where $1\le i\le s$ and $j\ge 1$. Notice that since $f\ne 1$, we have $\lambda_s\ne 0$. 
And by abuse of notation, we define $\max(f)\coloneqq \max\Set{\sigma(k):\lambda_k=\lambda_s}$. 

\begin{Lemma}
    \label{clm:G}
    The set $\calG_f$ is the disjoint union of the following three types of sets:
    \begin{itemize}
        \item $\{z_{\sigma(i),j}:1\le j < \lambda_i \text{ and } \lambda_i>0\}$,
        \item $\{z_{\sigma(i),j}:1\le j < \lambda_s \text{ and } \lambda_i=0\}$, and
        \item $\{z_{\sigma(i),\lambda_s}:\sigma(i)<\max(f) \text{ and } \lambda_i=0\}$.
    \end{itemize}
\end{Lemma}

A quick example illustrating the linear generating set $\calG_f$ above is given in \Cref{exam:the-one}.

\begin{proof}
    To confirm the assertion, roughly speaking, we have two cases.

    \begin{enumerate}[A]
        \item \label{Item-1}
            Firstly, we consider the case when $\lambda_i>0$.  We will show that $z_{\sigma(i),j} \in \calG_f$ if and only if $1\le j< \lambda_i$.
            \begin{enumerate}[i]
                \item  \label{Item-1-i} When $1\le j< \lambda_i$, we look at the squarefree monomial $f'\coloneqq f z_{\sigma(i),j} /z_{\sigma(i), \lambda_i}$. By (S-2), we have $f'\in J$. Thus, we can find suitable $f''\in G(J)$ that divides $f'$. Obviously, we have $\suppdeg(f'')\le \suppdeg(f')=\suppdeg(f)$ and $w(f'')\le w(f')<w(f)$. It follows that  $f''\prec f$.  Now, $z_{\sigma(i),j}f \in \braket{f'}\subseteq \braket{f''}\subseteq J_f$.  Consequently, $z_{\sigma(i),j} \in J_f: f$, as expected.
                \item  Suppose that $j= \lambda_i$.  It follows from \Cref{rmk:colon} that $z_{\sigma(i),\lambda_i}\notin \calG_f$. 
                \item \label{Item-1-iii} When $j> \lambda_i$, we suppose for contradiction that $z_{\sigma(i),j} \in \calG_f$. Since now $fz_{\sigma(i),j}\in J_f$, there exists some $f''\in G(J_f)$ such that $f''$ divides $fz_{\sigma(i),j}$. If $z_{\sigma(i),j}$ does not divide this $f''$, then $f''$ divides $f$ directly, contradicting the minimality of $f$. Thus, $z_{\sigma(i),j}$ divides $f''$. Since $f''$ is a sparse monomial, $z_{\sigma(i),\lambda_i}$ does not divide $f''$. Consequently, $f''$ divides the squarefree monomial $f'\coloneqq f z_{\sigma(i),j} /z_{\sigma(i), \lambda_i}$. Now, $f'\in J$ such that $\suppdeg(f')=\suppdeg(f)$ and $|\bdlambda(f')|>|\bdlambda(f)|$. Hence, we don't have $f'\prec f$. Consequently, $f''\ne f'$ and $\suppdeg(f'')<\suppdeg(f')=\suppdeg(f)$.  Next, we look at $f'''\coloneqq f''z_{\sigma(i), \lambda_i}/z_{\sigma(i),j}$, which belongs to $J$ by (S-2). Notice that $f'''$ properly divides $f$, as $f''$ properly divides $f'$. This contradicts the minimality of $f$. Therefore, $z_{\sigma(i),j} \notin J_f: f$, as expected.
            \end{enumerate}

        \item \label{Item-2}
            Next, we consider the case when $\lambda_i=0$. %Notice that $\lambda_s=\max(\bdlambda)>0$. 
            We will show that $z_{\sigma(i),j}\in \calG_f$ if and only if $1\le j < \lambda_s$, or $j=\lambda_s$ and $\sigma(i)<\max(f)$.
            %$z_{\sigma(i),\lambda_s}\in \calG_f$, when , and $z_{\sigma(i),\lambda_s}\notin \calG_f$ when $\sigma(i)>\max(f)$. For any fixed $j$, we show that  and $z_{\sigma(i),j} \notin \calG_f$ when $j>\lambda_s$. 

            \begin{enumerate}[i]
                \item When $1\le j < \lambda_s$, we look at the squarefree monomial $f'\coloneqq f z_{\sigma(i),j} /z_{\sigma(s), \lambda_s}$, which belongs to $J$ by \Cref{lem-2}. Then, we can argue as in the item \ref{Item-1-i} of \ref{Item-1} to see that $z_{\sigma(i),j}\in J_f:f$.
                    %Since $\suppdeg(f')=\suppdeg(f)$ and $w(f')<w(f)$, it follows that $f'\prec f$ and $f'\in J_f$. Consequently, we have $z_{\sigma(i),j}\in J_f:f$, as expected. 

                \item \label{Item-ii}
                    If $j>\lambda_s$, 
                    we suppose for contradiction that $z_{\sigma(i),j} \in \calG_f$.
                    As in the item \ref{Item-1-iii} of \ref{Item-1}, there exists some $f''\in G(J_f)$ such that $f''$ divides $fz_{\sigma(i),j}$, and $z_{\sigma(i),j}$ divides $f''$ by the minimality of $f$. Furthermore, since $f''\prec f$, we have $\suppdeg(f'')\le \suppdeg(f)$ and $f''$ properly divides $fz_{\sigma(i),j}$. We have two subcases to consider here.
                    \begin{enumerate}[a]
                        \item If $\suppdeg(f'')= \suppdeg(f)$, there is some $i'\ne i$ such that $\lambda_{i'}>0$ and $f''=fz_{\sigma(i),j}/z_{\sigma(i'), \lambda_{i'}}$. Since $\suppdeg(f'')=\suppdeg(f)$ while $w(f'')>w(f)$, this contradicts the condition $f''\prec f$.
                        \item \label{Item-ii-b}
                            If $\suppdeg(f'')<\suppdeg(f)$, there exists some subset 
                            \[
                                \calI \subseteq \Set{k:1\le k\le s, \lambda_k>0},
                            \]
                            such that $|\calI|\ge 2$, $i\notin \calI$ and $f''=f z_{\sigma(i),j}\big/ \prod_{k\in \calI} z_{\sigma(k),\lambda_k}$.  We can fix an arbitrary $k_0\in \calI$ and look at 
                            \[
                                f'''\coloneqq f'' z_{\sigma(k_0),\lambda_{k_0}}/z_{\sigma(i),j}=f\bigg/\prod_{k\in \calI\setminus \{k_0\}} z_{\sigma(k),\lambda_k}.
                            \]
                            Since $j>\lambda_{k_0}$, it follows from \Cref{lem-1} that $f'''\in J$ by considering the pair $(f'',f''')$. However, $f'''$ divides $f$ properly. This contradicts the minimality of $f$.
                    \end{enumerate}
                    Since we get contradiction in both subcases, $z_{\sigma(i),j}\notin \calG_f$, as expected.

                \item \label{Item-iii}
                    If $\sigma(i)<\max(f)$, we look at the squarefree monomial $f'\coloneqq f z_{\sigma(i),\lambda_s} /z_{\max(f), \lambda_s}$, which belongs to $J$ by (S-3). Notice that although $\bdlambda(f')=\bdlambda(f)$, we have $f'>_{\lex} f$. Thus, $f'\prec f$, which in turn implies  $z_{\sigma(i),\lambda_s}\in \calG_f$.

                \item It remains to consider the situation when $\sigma(i)\ge \max(f)$. Since $\sigma(i)=0<\lambda_s$, this amounts to saying that $\sigma(i)>\max(f)$. Now, we suppose for contradiction that $z_{\sigma(i),\lambda_s} \in \calG_f$.
                    Like in the case \ref{Item-ii}, we can find some $f''\in G(J_f)$ such that $z_{\sigma(i),\lambda_s}$ divides $f''$ and $f''$ properly divides $fz_{\sigma(i),\lambda_s}$.
                    \begin{enumerate}[a]
                        \item If $\suppdeg(f'')=\suppdeg(f)$, there is some $i'\ne i$ such that $\lambda_{i'}>0$ and $f''=fz_{\sigma(i),\lambda_s}/z_{\sigma(i'),\lambda_{i'}}$. We always have $\lambda_{i'}\le \lambda_s$. If $\lambda_{i'} <\lambda_s$, then $w(f'')>w(f)$, contradicting the assumption that $f''\prec f$.  On the other hand, if $\lambda_{i'}= \lambda_s $, then $\bdlambda(f'')=\bdlambda(f)$. As $\sigma(i)> \max(f)\ge \sigma(i')$, we have $f>_{\lex} f''$. And this again contradicts the assumption that $f''\prec f$.

                        \item If $\suppdeg(f'')<\suppdeg(f)$, the proof will be similar to that in item \ref{Item-ii-b} of \ref{Item-ii}, with $j$ there replaced by $\lambda_s$ here. Note that we can still apply \Cref{lem-1}.
                    \end{enumerate}
                    Since we get contradiction in both subcases, $z_{\sigma(i),\lambda_s}\notin \calG_f$, as expected.
            \end{enumerate}

    \end{enumerate}

    Thus, we have proved the \Cref{clm:G}.
    %Furthermore, if $J_f: f$ is linear, it must be generated by $\calG_f$.
\end{proof}

%\subsection{Linearity of the colon ideal}
%In this subsection, we will 
By the definition of $\calG_f$, we have $\braket{\calG_f}\subseteq J_f:f$.
The remaining task is to show this is actually equality.

\begin{Lemma}
    \label{lem:containment}
    We have $\braket{\calG_f}=J_f:f$ in general. 
\end{Lemma}

\begin{proof}
    Suppose for contradiction that we only have $\braket{\calG_f}\subsetneq
    J_f:f$.  Then, we can find suitable  $g\in G(J_f:f)$ such that $g\notin
    \braket{\calG_f}$. We have mentioned in \Cref{rmk:colon} that this $g$ is a
    sparse squarefree monomial which is coprime to $f$. Thus, there exist some
    nonempty $\calJ\subseteq [s]\coloneqq \{1,2,\dots,s\}$ and positive
    integers $\lambda_j'$ for each $j\in \calJ$, such that $g=\prod_{j\in
        \calJ}z_{\sigma(j),\lambda_j'}$.  Thus, as $g\notin \calG_f$ and $g$ is
    coprime to $f$, for each $j\in \calJ$, we have:
    \begin{enumerate}[i]
        \item if $\lambda_j>0$, then $\lambda_j'>\lambda_j$;
        \item if $\lambda_j=0$, then either $\lambda_j'>\lambda_s$, or $\lambda_j'=\lambda_s$ and $\sigma(j)>\max(f)$.
    \end{enumerate}
    As $gf\in J_f$, we can find some $f'\in G(J_f)$ dividing $gf$. Since $g$ is a minimal generator of $J_f:f$ while $f'$ is a sparse squarefree monomial, we must have $f'= g \prod_{j\in \calJ'} z_{\sigma(j),\lambda_j}$
    for some $\calJ'\subset \Set{j\in [s]\setminus \calJ:\lambda_j>0}$.
    Furthermore, as $f'\prec f$, we have $\suppdeg(f')\le \suppdeg(f)$. Now, we have two cases to discern.
    \begin{enumerate}[a]
        \item The first case is when $\suppdeg(f')<\suppdeg(f)$, which amounts to saying that 
            \[
                \left| \Set{j\in \calJ:\lambda_j=0} \right| < \left| \Set{j\notin \calJ\sqcup \calJ':\lambda_j>0} \right|.
            \]
            Here, $\sqcup$ indicates that the union is disjoint. We may take an arbitrary subset $\calJ''\subset \Set{j\notin \calJ\sqcup \calJ':\lambda_j>0} $ with a bijection $\varphi$ from $\calJ''$ to  $\Set{j\in \calJ:\lambda_j=0}$. Since $f'\in J$, the product element
            \[
                \prod_{j\in \calJ:\lambda_j>0}z_{\sigma(j),\lambda_j'} \cdot \prod_{j\in \calJ'}z_{\sigma(j),\lambda_j} \cdot \prod_{j\in \calJ''} z_{\sigma(j),\lambda_{\varphi(j)}'} \in J
            \]
            by the symmetry condition (S-3). 
            In the first part of this product, each $\lambda_j'>\lambda_j$. Meanwhile, in the third part, each $\lambda_{\varphi(j)}'\ge \lambda_s\ge \lambda_j$. Thus, by the shiftiness condition (S-2), we have 
            \[
                f''\coloneqq \prod_{j\in \calJ:\lambda_j>0}z_{\sigma(j),\lambda_j} \cdot \prod_{j\in \calJ'}z_{\sigma(j),\lambda_j} \cdot \prod_{j\in \calJ''} z_{\sigma(j),\lambda_{j}}\in J.
            \]
            As $\suppdeg(f'')=\suppdeg(f')<\suppdeg(f)$, $f''$ properly divides $f$,  contradicting the minimality of $f$.
        \item The second case is when $\suppdeg(f')=\suppdeg(f)$. 
            If we apply the notation in the previous case, then this simply means $\Set{j\in [s]:\lambda_j>0}$ is a subset of the disjoint union $\calJ\sqcup \calJ'\sqcup \calJ''$.  It then follows from $f'\prec f$ that $w(f')\le w(f)$, which is just 
            \[
                \sum_{\substack{j\in \calJ\\ \lambda_j>0}}\lambda_j' + \sum_{j\in \calJ'}\lambda_j+\sum_{j\in \calJ''} \lambda_{\varphi(j)}'
                \le
                \sum_{\substack{j\in \calJ\\ \lambda_j>0}}\lambda_j + \sum_{j\in \calJ'}\lambda_j+\sum_{j\in \calJ''} \lambda_{j}.
            \]
            By the comparisons stated in the previous case, this is possible
            only when all of the following hold:
            \begin{enumerate}[1]
                \item \label{J-item-1} there is no $j\in \calJ$ such that
                    $\lambda_j>0$; 
                \item \label{J-item-2} for each $j\in \calJ''$,
                    $\lambda_j=\lambda_s=\lambda_{\varphi(j)}'$.
            \end{enumerate}
            Consequently, $\bdlambda(f')=\bdlambda(f)$ and
            $\varphi(\calJ'')=\Set{j\in \calJ:\lambda_j=0}=\calJ$. Furthermore,
            since
            \[
                \calJ'\sqcup \calJ'' \subseteq \Set{j:\lambda_j>0} \subseteq \calJ\sqcup \calJ'\sqcup \calJ'',
            \]
            it follows again from the above item \ref{J-item-1} that
            \[
                %\left| \calJ \right|= \left| \calJ'' \right|= \left| \Set{j\in [s]\setminus\calJ':\lambda_j>0} \right|.
                 \Set{j\in [s]:\lambda_j>0} =\calJ' \sqcup \calJ''.
            \]
            And as $\calJ=\varphi(\calJ'')$, we have $\lambda_j'=\lambda_s$ for
            each $j\in \calJ$, by the above item \ref{J-item-2}.  Hence, we can
            observe that
            \begin{itemize}
                \item for each $j\in [s]\setminus \calJ'$ with $\lambda_j\ne
                    \lambda_j'$, $j$ belongs to $\calJ''$ or $\calJ$, which
                    means that $(\lambda_j,\lambda_j')$ equals either
                    $(\lambda_s,0)$ or $(0,\lambda_s)$ correspondingly;
                \item for each $j\in [s]\setminus \calJ'$ with
                    $(\lambda_j,\lambda_j')=(0,\lambda_s)$, one has $j\in
                    \calJ$ and $\sigma(j)>\max(f)$ by the property of $g$.
            \end{itemize}
            Since $\calJ\ne\varnothing$, it follows that $f>_{\lex}f'$,
            contradicting to the fact $f'\in G(J_f)$.
    \end{enumerate}
    Therefore, we have shown the claimed equality $\braket{\calG_f}=J_f:f$. 
\end{proof}

Now, $J_f:f$ is generated by the variables in $\calG_f$. In particular, it is
linear. Since this holds for each $f\in G(J)$, $J$ has linear quotients, as
expected.

\subsection{Betti Numbers}
In the following, $J$ is a sparse symmetric shifted ideal. Like \cite[Corollary
5.7]{MR4108337}, we want to calculate the graded Betti number of $J$. The main
arguments for that purpose follow closely from those in \cite[Section
5]{MR4108337}. Some necessary details are included here, mainly to keep the
exposition self-contained.

Let $\Lambda(J)\coloneqq \Set{\bdlambda(f):f\in G(J)}$ be the set of partitions
of the minimal monomial generators of $J$. Without loss of generality, we may
assume that $\Lambda(J)=\{\bdlambda^{(1)},\dots ,\bdlambda^{(t)}\}$ with
$\bdlambda^{(1)}\prec \cdots \prec \bdlambda^{(t)}$.  Here, $\prec$ is the
total order we introduced in the proof of \Cref{SSSisLQ}.  For each $k$ with
$1\le k\le t$, let $J_{\le k}\subseteq J$ be the sparse symmetric shifted ideal
with $\Lambda(J_{\le k})=\{\bdlambda^{(1)}, \dots ,\bdlambda^{(k)}\}$.  Thus,
we have a short exact sequence
\[
    0\to J_{\le k-1}\to J_{\le k}\to J_{\le k}/J_{\le k-1} \to 0.
\]
Notice that the quotient in the above sequence can be better recognized as
\[
    J_{\le k}/J_{\le k-1}\cong \braket{\sigma(\bdz^{\bdlambda^{(k)}}):\sigma\in\frakS_s} \left/ \left(\braket{\sigma(\bdz^{\bdlambda^{(k)}}):\sigma\in\frakS_s} \cap J_{\le k-1}\right).\right.
\]

\begin{Lemma}
    %\label{betti-paper-lem-5.6}
    Assume that $\bdlambda^{(k)}=(\lambda_1,\dots,\lambda_s)$. Then, the
    intersection ideal $K\coloneqq\braket{\frakS_s\bdz^{\bdlambda^{(k)}}} \cap
    J_{\le k-1}$ is the sum of the ideal
    \[
        K_1\coloneqq \braket{\sigma(\bdz^{\bdlambda^{(k)}}z_{i,j}): \sigma\in\frakS_s, \lambda_i>0\text{ and } 1\le j<\lambda_i }
    \]
    with the ideal
    \[
        K_2\coloneqq \braket{\sigma(\bdz^{\bdlambda^{(k)}}z_{i,j}): \sigma\in\frakS_s, \lambda_i=0\text{ and } 1\le j<\lambda_s  }.
    \]
\end{Lemma}

\begin{proof}
    We first show the containment ``$K \supseteq K_1+K_2$''. Let
    $f=\bdz^{\bdmu}$ be the first monomial in $G(J_{\le k})\setminus G(J_{\le
        k-1})$ with respect to the order $\prec$ on $G(J)$. Then, for some
    $\tau\in \frakS_s$, we have $f=\tau(\bdz^{\bdlambda^{(k)}})=\prod_i
    z_{\tau(i),\lambda_i}$, and $J_f:f=J_{\le k-1}:f=\braket{\calG_f}$.  Notice
    that in this case,
    \begin{equation}
        \calG_f=
        \Set{z_{\tau(i),j}: \lambda_i>0\text{ and } 1\le j < \lambda_i}\cup \Set{z_{\tau(i),j}: \lambda_i=0\text{ and } 1\le j < \lambda_s}.
        %\tag{$\ast$}
        \label{eqn:short-G_f}
    \end{equation}
    Clearly, the lexicographic order is used here to make the third type of sets in
    \Cref{clm:G} empty, via the item \ref{Item-2}\ref{Item-iii} in its proof.
    Now, for each $z_{\tau(i),\lambda_j}\in \calG_f$, one has $f
    z_{\tau(i),\lambda_j}=\tau(\bdz^{\bdlambda^{(k)}} z_{i,j})\in J_{\le k-1}$.
    And the containment ``$\supseteq$'' follows from the $\frakS_s$ actions.

    Next, we show the opposite containment ``$K \subseteq K_1+K_2$''. For this
    purpose, we take arbitrary squarefree monomial $g$ in
    $\braket{\frakS_s\bdz^{\bdlambda^{(k)}}} \cap J_{\le k-1}$. We will also
    utilize the generator $f=\tau(\bdz^{\bdlambda^{(k)}})$ in the previous
    paragraph. Now, $\tau'(g)\in \braket{f}$ for some $\tau'\in \frakS_s$.
    Thus, $\tau'(g)=g_1f$ for some squarefree monomial $g_1$. Since
    $\tau'(g)\in J_{\le k-1}$, $g_1\in J_{\le k-1}:f=\braket{\calG_f}$, which
    was described in the above equation \eqref{eqn:short-G_f}. Now the
    containment ``$\subseteq$'' follows from the $\frakS_s$ actions.
\end{proof}

In the following, for each partition $\bdlambda=(\lambda_1,\dots,\lambda_s)$,
we will write $N^{\bdlambda}$ for the quotient module
\begin{align*}
    \frac{\braket{\sigma(\bdz^{\bdlambda}): \sigma\in\frakS_s}}{
        \braket{\sigma(\bdz^{\bdlambda} z_{i,j}): \lambda_i>0,1\le j<\lambda_i, \sigma\in\frakS_s}+\braket{\sigma(\bdz^{\bdlambda} z_{i,j}): \lambda_i=0,1\le j<\lambda_s, \sigma\in\frakS_s}}.
\end{align*}
%Thus, \Cref{betti-paper-lem-5.6} says that 
\begin{Corollary}
    One has $J_{\le k}/J_{\le k-1}\cong N^{\bdlambda^{(k)}}$.
\end{Corollary}

Next, we want to describe the $\KK$-linear structure of $N^{\bdlambda}$. Let
$N_{num}^{\bdlambda}$ and $N_{denom}^{\bdlambda}$ be the squarefree monomial
ideals in the numerator and denominator of $N^{\bdlambda}$ respectively.  Since
$N^{\bdlambda}$ is a squarefree module, it suffices to take an arbitrary squarefree
monomial $g\in N_{num}^{\bdlambda}$ and ask when $g\notin
N_{denom}^{\bdlambda}$. For each $i\in [s]$, if some $z_{i,j}$
divides $g$, let $\mu_i=\min\Set{j: \text{$z_{i,j}$ divides $g$}}$. If no such
$z_{i,j}$ exists for this chosen $i$, let $\mu_i=0$. Thus, we have a tuple
$\bdmu=(\mu_1,\dots,\mu_s)\in \ZZ_{\ge 0}^s$ for this squarefree monomial $g$.
Additionally, we may assume that the partition associated with $\bdz^{\bdmu}$
is $\bdlambda'=(\lambda_1',\dots,\lambda_s')$.  Suppose that
$\theta=\suppdeg(\bdlambda)$ and $\theta'=\suppdeg(\bdlambda')$. Thus,
$\lambda_{s-\theta}=0<\lambda_{s-\theta+1}$ and similarly
$\lambda_{s-\theta'}'=0<\lambda_{s-\theta'+1}'$.

\begin{Lemma}
    \label{lem:N-Lambda}
    With the notations as above, the squarefree monomial $g$ belongs to
    $N_{num}^{\bdlambda}\setminus N_{denom}^{\bdlambda}$ if and only if 
    \begin{equation}
        \text{$\theta\le \theta'\quad$ and $\quad\lambda_{s-\theta+k}=\lambda_{s-\theta'+k}'$ for $k=1,2,\dots,\theta$.}
        %\tag{\dag}
        \label{eqn:dag}
    \end{equation} 
\end{Lemma}

\begin{proof}
    Suppose that the squarefree monomial $g$ belongs to
    $N_{num}^{\bdlambda}\setminus N_{denom}^{\bdlambda}$. By symmetry, we may
    assume that $\bdz^{\bdlambda}$ divides $g$. We observe the following two
    facts.
    \begin{enumerate}[a]
        \item Firstly, whenever $\lambda_i>0$, $z_{i,\lambda_i}$ divides $g$.
            Therefore, $\mu_i \le \lambda_i$. However, when $\mu_i< \lambda_i$
            for one such $i$, then $g\in N_{denom}^{\bdlambda}$ as well. Thus,
            we have precisely $\mu_i=\lambda_i$ whenever $\lambda_i>0$. 
        \item Secondly, for each $i$ with $\lambda_i=0$, if $1\le
            \mu_i<\lambda_s$, then $g\in N_{denom}^{\bdlambda}$ as well. Thus
            either $\mu_i=0$ or $\mu_i\ge \lambda_s$. 
    \end{enumerate}
    From these observations, we see immediately that condition \eqref{eqn:dag} holds.

    Conversely, suppose that condition \eqref{eqn:dag} is satisfied. It is
    clear that $\bdz^{\bdlambda}$ divides $\tau(\bdz^{\bdlambda'})$ for some
    permutation $\tau\in\frakS_s$. Meanwhile,
    $\tau'(\bdz^{\bdlambda'})=\bdz^{\bdmu}$ divides $g$ for some permutation
    $\tau'\in \frakS_s$.  Thus, $g\in \braket{\sigma(\bdz^{\bdlambda}):
        \sigma\in\frakS_s}=N_{num}^{\bdlambda}$.  As for testing whether
    $g\in N_{denom}^{\bdlambda}$, we are reduced to the following two cases by symmetry.
    \begin{enumerate}[i]
        \item Suppose that $g$ is divisible by $\bdz^{\bdlambda} z_{i_0,j}$ for
            some $1\le j<\lambda_{i_0}$ with $\lambda_{i_0}>0$. Then, whenever
            $\lambda_{i}>0$, we have $\mu_{i}\le \lambda_{i}$. Furthermore,
            $\mu_{i_0}<\lambda_{i_0}$. This will break the assumptions in
            \eqref{eqn:dag}.
        \item Suppose that $g$ is divisible by $\bdz^{\bdlambda} z_{i_0,j}$ for
            some $1\le j<\lambda_s$ with $\lambda_{i_0}=0$. Then, whenever
            $\lambda_{i}>0$, we have $\mu_{i}\le \lambda_{i}$. Furthermore,
            $1\le \mu_{i_0}<\lambda_{s}$. This still breaks the assumptions in
            \eqref{eqn:dag}.
    \end{enumerate}
    Thus, $g\in N_{num}^{\bdlambda}\setminus N_{denom}^{\bdlambda}$.
    %    Conversely, suppose that the condition (\dag) is satisfied. It is clear that
    %    $\bdz^{\bdlambda}$ divides $\bdz^{\bdlambda'}$ while $\bdz^{\bdlambda'}$ divides $g$. Thus, $g\in \braket{\sigma(\bdz^{\bdlambda}): \sigma\in\frakS_s}=N_{num}^{\bdlambda}$. It remains to show that $g\notin N_{denom}^{\bdlambda}$. Suppose for contradiction that this is not true. Then, by symmetry, we only need to consider the following two cases.
    %    \begin{enumerate}[a]
    %       \item Suppose that $g$ is divisible by $\bdz^{\bdlambda} z_{i_0,j}$ for some $1\le j<\lambda_{i_0}$ with $\lambda_{i_0}>0$. Then, whenever $\lambda_{i}>0$, we have $\mu_{i}\le \lambda_{i}$. Furthermore, $\mu_{i_0}<\lambda_{i_0}$. This will break the assumptions in (\dag).
    %        \item Suppose that $g$ is divisible by $\bdz^{\bdlambda} z_{i_0,j}$ for some $1\le j<\lambda_s$ with $\lambda_{i_0}=0$. Then, whenever $\lambda_{i}>0$, we have $\mu_{i}\le \lambda_{i}$. Furthermore, $1\le \mu_{i_0}<\lambda_{s}$. This still breaks the assumptions in (\dag).
    %   \end{enumerate}
    %   Thus, $g\in N_{num}^{\bdlambda}\setminus N_{denom}^{\bdlambda}$.
\end{proof}

It is time to introduce a series of new notations.
\begin{Notation}
    Given any partition $\bdlambda=(\lambda_1,\dots,\lambda_s)$, we may assume that 
    \[
        \bdlambda=(0,\dots,0,\lambda_{k_1}, \lambda_{k_1+1},\dots,\lambda_{k_2},\lambda_s,\dots,\lambda_s) 
    \]
    where $\lambda_{k_1}\ne 0$ and $\lambda_{k_2}<\lambda_s$.  
    \begin{itemize}
        \item Let $\pi_{\bdlambda}\coloneqq
            (\lambda_{k_1},\lambda_{k_1+1},\dots,\lambda_{k_2})$ be the
            non-zero and non-maximal piece of $\bdlambda$. Write
            $p(\bdlambda)\coloneqq k_2-k_1+1=\left|\Set{j:1\le
                    \lambda_j<\lambda_s}\right|$ for the length of this piece
            and $r(\bdlambda)\coloneqq\left|\Set{j:\lambda_j=\lambda_s}\right|$
            for the length of the remaining maximal piece. 
        \item  For any $p(\bdlambda)$-subset
            $A=\{i_1,i_2,\dots,i_{p(\bdlambda)}\}\subseteq [s]$ with
            $i_1<i_2<\cdots <i_{p(\bdlambda)}$, and any positive integer $a$,
            let $Z_{A,a}\coloneqq \KK[z_{i,j}\in Z:i\in A \text{ and } j\ge
            a]$. When $a=1$, we may write this subring simply as $Z_A$ and write
            $\bdz_{A}^{\pi_{\bdlambda}}\coloneqq \prod_{j=1}^{p(\bdlambda)}
            z_{i_{j},\lambda_{k_1+j-1}}\in Z_A$.  
        \item  With respect to the subset $A$ above, we write
            ${A}^{\complement}\coloneqq [s]\setminus A$ for the complement set,
            and similarly introduce $Z_{{A}^\complement,a}\coloneqq
            \KK[z_{i,j}\in Z:i\in {A}^{\complement} \text{ and } j\ge a]$ for
            the positive integer $a$.
        \item For the complement set $A^{\complement}$, we will write
            $V_{A^{\complement},\bdlambda}$ for the monomial ideal in
            $Z_{A^{\complement},\lambda_s}$ generated by the sparse squarefree
            monomials of the form $z_{t_1,\lambda_s}z_{t_2,\lambda_s}\cdots
            z_{t_{r(\bdlambda)},\lambda_s}$ with
            $t_1,t_2,\dots,t_{r(\bdlambda)}$ distinct in $A^{\complement}$. 
    \end{itemize}
\end{Notation}

Notice that we have a natural group action of $\frakS_{A}$ on $Z_{A,a}$, which
can be derived from the action of $\frakS_s$ on $Z_{A,a}$.

\begin{Remark}
    Fix a sparse squarefree monomial $m= \prod_{j=1}^{p(\bdlambda)}
    z_{i_{j},\mu_{j}} \in \frakS_A \bdz_{A}^{\pi_{\bdlambda}}$ with
    $\bdmu=(\mu_1,\dots,\mu_{p(\bdlambda)})$.  Surely we have $mm'\in \frakS_s
    \bdz^{\bdlambda}$ for each sparse monomial $m'\in
    G(V_{A^{\complement},\bdlambda})$.  It is not difficult to deduce from
    \Cref{lem:N-Lambda} that the annihilating ideal of the image of $mm'$ is
    given by
    \begin{equation}
        \braket{z_{i,j}:
            \text{either $i\in A$ and $1 \le j<\mu_i$, or $i\notin A$ and $1\le j<\lambda_s$}
        },
        %\tag{\ddag}
        \label{eqn:m-ann}
    \end{equation} 
    which is independent of the choice of $m'$.
\end{Remark}

\begin{Notation}
Inspired by the description of \eqref{eqn:m-ann}, we introduce the following
subrings, submodule, and ideals. 
    \begin{itemize}
        \item For any given $m= \prod_{j=1}^{p(\bdlambda)} z_{i_{j},\mu_{j}}
            \in \frakS_A \bdz_{A}^{\pi_{\bdlambda}}$ with
            $\bdmu=(\mu_1,\dots,\mu_{p(\bdlambda)})$, let 
            \[
                Z_m^{\vartriangle}\coloneqq \KK[z_{i,j}:
                \text{either $i\in A$ and $1 \le j<\mu_i$, or $i\notin A$ and $1\le j<\lambda_s$}],
            \]
            and write $\frakm_{Z_m^{\vartriangle}}$ for its graded maximal
            ideal.  Similarly, let 
            \[
                Z_m^\triangledown \coloneqq \KK[z_{i,j}:
                \mu_i \le j\le b \text{ and } i\in A]\subset Z_A.
            \]
            Roughly speaking, the variables in $Z_m^{\vartriangle}$ are in the
            upper half of the canonical matrix $\bdZ$ consider in \Cref{rmk:sss},
            while the variables in $Z_m^{\triangledown}$ are in the lower half
            of the matrix.
        \item For this given $m$, we write $N_{A,m}^{\bdlambda}$ for the
            submodule generated by the image of
            $mV_{A^{\complement},\bdlambda}$ in $N^{\bdlambda}$.
        \item Additionally, we write $V_{A^{\complement},\bdlambda}^{*}$ for
            the squarefree Veronese ideal of degree $r(\bdlambda)$ in
            $\KK[z_{i,\lambda_s}:i\in A^{\complement}]$, which is generated by
            the squarefree monomials of degree $r(\bdlambda)$ in this ring.
    \end{itemize}
\end{Notation}

\begin{Remark}
    It is clear that 
    \begin{equation}
        V_{A^{\complement},\bdlambda}=V_{A^{\complement},\bdlambda}^{*}Z_{A^{\complement},\lambda_s}\cong V_{A^{\complement},\bdlambda}^{*}\otimes_{\KK}Z_{A^{\complement},\lambda_s+1} 
        \label{eqn:Veronese}
    \end{equation}
    is an extended ideal from $V_{A^{\complement},\bdlambda}^{*}$.  Meanwhile,
    we also have that multigraded isomorphism
    \begin{equation}
        N_{A,m}^{\bdlambda}\cong Z_m^{\vartriangle}/\frakm_{Z_m^{\vartriangle}} \otimes_\KK mZ_m^\triangledown\otimes_{\KK} V_{A^{\complement},\bdlambda}^{*}\otimes_\KK Z_{A^{\complement},\lambda_s+1}.
        \label{eqn:Naml}
    \end{equation}
\end{Remark}

Thus, we can deduce from \Cref{lem:N-Lambda} and the multidegree information in
\eqref{eqn:Naml} the following conclusion.

\begin{Corollary}
    We have the following decomposition as multigraded $Z$-modules:
    \[
        N^{\bdlambda}=\bigoplus_{\substack{A\subset[s]\\ |A|=p(\bdlambda)}}\bigoplus_{m\in \frakS_{A}\bdz_A^{\pi_{\bdlambda}}} 
        N_{A,m}^{\bdlambda}.
        %\overline{m} Z_m V_{A^{\complement},\bdlambda}.
    \]
\end{Corollary}

Notice that the monomial $m$ sits only in $Z_A$ and has degree
$\deg(\bdz_A^{\pi_{\bdlambda}})=p(\bdlambda)$.  Since
$V_{A^{\complement},\bdlambda}^{*}$ is the squarefree Veronese ideal of degree
$r(\bdlambda)$ in a suitable ring. It follows from \eqref{eqn:Veronese} that
$V_{A^{\complement},\bdlambda}$ in $Z_{A^{\complement},\lambda_s}$ has a linear resolution with regularity
$r(\bdlambda)$.  Consequently, both $N_{A,m}^{\bdlambda}$ and $N^{\bdlambda}$
have linear resolutions with regularity
$p(\bdlambda)+r(\bdlambda)=\suppdeg(\bdlambda)$.

Now, we can state the second main result of this section.

\begin{Theorem}
    If $J$ is a sparse symmetric shifted ideal, then as $\KK[\frakS_s]$-modules, we have
    \[
        \Tor_i(J,\KK)_{i+d}\cong \bigoplus_{\substack{\bdlambda\in \Lambda(J)\\ \suppdeg(\bdlambda)=d}}\Tor_i(N^{\bdlambda},\KK).
    \]
\end{Theorem}

Its proof is literally the same as the argument for \cite[Theorem 5.5]{MR4108337}.  There is no need repeating it here. 

It remains to give a closed formula of the graded Betti numbers of $J$ using the information from $\Lambda(J)$. To achieve this, consider the case when $\suppdeg(\bdlambda)=d$, $A\subset[s]$ with $|A|=p(\bdlambda)$, and $m=\prod_{j=1}^{p(\bdlambda)} z_{i_{j},\mu_{j}} \in \frakS_{A}\bdz_A^{\pi_{\bdlambda}}$.
%Inspired by the equality \eqref{eqn:m-ann} for this monomial $m$, 
Notice that the dimension of $Z_m^{\vartriangle}$ is precisely
\begin{align*}
    p^\vartriangle(\bdlambda)&\coloneqq |\pi_{\bdlambda}|-p(\bdlambda)+(\lambda_s-1)(s-p(\bdlambda))\\
    %&=|\bdlambda|-\suppdeg(\bdlambda)+(\lambda_s-1)(s-\suppdeg(\bdlambda))\\
    &=|\bdlambda|+s(\lambda_s-1)-\lambda_s\suppdeg(\bdlambda).
\end{align*}
Thus, we can deduce from \eqref{eqn:Naml} that 
\begin{align*}
    \beta_{i,i+d}(N_{A,m}^{\bdlambda})&=\beta_i(Z_m^{\vartriangle}/\frakm_{Z_m^{\vartriangle}} \otimes_\KK Z_m^\triangledown\otimes_{\KK} V_{A^{\complement},\bdlambda}^{*}\otimes_\KK Z_{A^{\complement},\lambda_s+1})\\
    &=\sum_{k+\ell=i} \beta_\ell(Z_m^{\vartriangle}/\frakm_{Z_m^{\vartriangle}})\beta_k (V_{A^{\complement},\bdlambda}^{*})\\
    &=\sum_{k+\ell=i}\binom{p^{\vartriangle}(\bdlambda)}{\ell}\binom{s-p(\bdlambda)}{r(\bdlambda)+k}\binom{r(\bdlambda)+k-1}{k};
\end{align*}
see also \cite[Theorem 2.1]{MR4105544}. 

Meanwhile, we will define the term $\type(\pi_{\bdlambda})\coloneqq (t_1,\dots,t_{\lambda_{s-1}})$ where $t_k=\left| \{j:\lambda_j=k\} \right|$ for $1\le k\le \lambda_{s-1}$. And it is clear that $\left| \frakS_A\bdz_{A}^{\pi_{\bdlambda}} \right|=\frac{p(\bdlambda)!}{\type(\pi_{\bdlambda})!}$, where $\type(\pi_{\bdlambda})!\coloneqq t_1!t_2!\cdots t_{\lambda_{s-1}}!$. Consequently, we have the following closed formula for Betti numbers.

\begin{Corollary}\label{SSSBetti}
    If $J$ is a sparse symmetric shifted ideal, then $\beta_{i,i+d}(J)$ is given by
    \[
        \sum_{\substack{\bdlambda\in\Lambda(J)\\\suppdeg(\bdlambda)=d}} \left( 
            \sum_{k+\ell=i}\frac{p(\bdlambda)!}{\type(\pi_{\bdlambda})!}\binom{p^{\vartriangle}(\bdlambda)}{\ell}\binom{s}{p(\bdlambda)}\binom{s-p(\bdlambda)}{r(\bdlambda)+k}\binom{r(\bdlambda)+k-1}{k}
        \right).
    \]
\end{Corollary}

A quick example using this formula is given in \Cref{exam:the-one}.

\section{Sequentially Cohen--Macaulay property of uniform $a$-fold product ideals}

In this section, we will show that the symbolic powers of the uniform $a$-fold product ideal are sequentially Cohen--Macaulay under some generic conditions. 
From now on, $\calF_{\bdz}$ will be the list of variables $\{z_1,\dots,z_s\}$ in the polynomial ring $T=\KK[z_1,\dots,z_s]$. And the field $\KK$ is infinite.

\subsection{Standard monomial case}
\begin{Proposition}
    \label{prop:mono-SCM-general}
    Suppose that $I=\bigcap_{k=1}^{t}I_{c_k,\calF_\bdz}^{(m_k)}\subset T$ is an
    intersection of symbolic powers of monomial star configuration ideals.
    Then, $I$ is sequentially Cohen--Macaulay.
\end{Proposition}

\begin{proof}
    To show that this ideal is sequentially Cohen--Macaulay, we will apply the
    polarization technique. And here is a list of facts that we will utilize.
    They are in Propositions 2.3, 2.5 and 4.11 of \cite{MR2184792}
    respectively.
    \begin{itemize}
        \item We have $\Pol(I_1\cap I_2)=\Pol(I_1)\cap \Pol(I_2) $.
        \item The polarization of $\braket{z_{i_1},\dots,z_{i_r}}^m$ has the
            following irredundant irreducible primary decomposition:
            \[
                \Pol(\braket{z_{i_1},\dots,z_{i_r}}^m)=\bigcap_{\substack{1\le e_j\le m,\\ \sum_j\! e_j\le m+r-1}} \braket{z_{i_1,e_1},\dots,z_{i_r,e_r}}.
            \]
        \item The ideal $I$ is sequentially Cohen--Macaulay if and only if
            $\Pol(I)$ is so.
    \end{itemize}
    %\begin{Lemma}
    %    [{\cite[Corollary 2.6]{MR2184792}}]
    %    Let $I$ be a monomial ideal in a polynomial ring $R=\KK[x_1,\dots,x_n]$ and let $\calP(I)$ be its polarization in $S=\KK[x_{i,j}]$. Then $(x_{i_1},\dots,x_{i_r})\in \Ass_R(R/I)$ if and only if $(x_{i_1,c_1},\dots,x_{i_r,c_r})\in \Ass_S(S/\calP(I))$ for some positive integers $c_1,\dots,c_r$. Moreover, if $(x_{i_1,c_1},\dots,x_{i_r,c_r})\in \Ass_S(S/\calP(I))$, then $(x_{i_1,b_1},\dots,x_{i_r,b_r})\in \Ass_S(S/\calP(I))$ for all $b_j$ such that $1\le b_j\le c_j$.
    %\end{Lemma}
    Now, we will apply the polarization to the ideals in the intersection 
    \[
        I=\bigcap_{k=1}^{t}I_{c_k,\calF_\bdz}^{(m_k)}=\bigcap_{k=1}^{t} \bigcap_{1\le i_1<\cdots<i_{c_k}\le s}\braket{z_{i_1},\dots,z_{i_{c_k}}}^{m_k}. 
    \]

    Notice that the squarefree monomial ideal $\Pol(I)$ is sequentially Cohen--Macaulay, if and only if the Alexander dual ideal  $\Pol(I)^{\vee}$ is componentwise linear, by \cite[Theorem 8.2.20]{MR2724673}. And now it suffices to show that $\Pol(I)^{\vee}$ has linear quotients, by \cite[Theorem 8.2.15]{MR2724673}.

    %Regarding the Alexander dual ideal $(I_a(\calF^b))^{\vee}$, I \alert{suspect} that it is symmetric shifted. If this is the case, then we are done by the following.
    %
    %\begin{Lemma}
    %    [{\cite[Theorem 3.2]{MR4108337}}]
    %    Symmetric shifted ideals have linear quotients.
    %\end{Lemma}
    Observe that $\Pol(I)^{\vee}$ is generated by monomials of the form $\prod_{j=1}^{c_k}z_{i_{j},e_{j}}$,
    where $1\le k\le t$, $1\le i_1<\cdots<i_{c_k}\le s$, $1\le e_j\le m_k$ and $\sum_j e_j\le m_k+c_k-1$. 
    It is clear that $\Pol(I)^{\vee}$ is sparse symmetric shifted. Thus, by \Cref{SSSisLQ}, this ideal has linear quotients, as expected. And this completes our proof for \Cref{prop:mono-SCM}.
\end{proof}

In the remaining of this paper, we will stay with the following setup.

\begin{Setting}
    \label{setting-monomial-case}
    Let $a,b$ be positive integers such that $(b-1)s+1<a\le bs$. And we introduce the positive integers 
    \[
        h\coloneqq bs-a+1, \quad
        c_0\coloneqq s-\floor{\frac{a-1}{b}} \quad 
        \text{and} \quad
        \mu_a^0\coloneqq a-b(s-c_0).  % s=n+1
    \]
\end{Setting}

Notice that if $a\le (b-1)s+1$, then $I_a(\calF_{\bdz}^b)$ is not saturated by \cite[Proposition 4.3]{arXiv:1912.04448}.
And under \Cref{setting-monomial-case}, the big height of $I_a(\calF_\bdz^b)$ is precisely $h$
by \cite[Remark 4.5]{arXiv:1912.04448}.
The main result of this subsection is the following.

\begin{Proposition}
    \label{prop:mono-SCM}
    Assume the notations in \Cref{setting-monomial-case}. Then, the symbolic powers of $I_a(\calF_\bdz^b)=I_a(z_1^b\cdots z_s^b)$ are all sequentially Cohen--Macaulay.
\end{Proposition}

\begin{proof}
    It follows from \cite[Theorem 4.6]{arXiv:1912.04448} that the symbolic power
    \begin{equation}
        I_a(\calF_\bdz^b)^{(m)}
        = \bigcap_{c=c_0}^{h} I_{c,\calF_\bdz}^{(m(\mu_a^0+b(c-c_0)))} 
        \label{eqn:IaDecomp}
    \end{equation}
    for each positive integer $m$. Now, it remains to apply \Cref{prop:mono-SCM-general}.
\end{proof}

\begin{Example}
    \label{exam:the-one} 
    Here, we consider a small example where $I=I_6(z_1^2z_2^2z_3^2z_4^2)$ in
    $\KK[z_1,z_2,z_3,z_4]$. Now, the squarefree monomial ideal
    $J=\Pol(I)^{\vee}$ in
    $\KK[z_{1,1},z_{2,1},z_{3,1},z_{4,1},z_{1,2},z_{2,2},z_{3,2},z_{4,2}]$ is
    sparse symmetric shifted, and has $22$ minimal monomial generators. We can
    arrange these monomials using the total order $\prec$ in the proof of
    \Cref{SSSisLQ}, and denote them by $u_1,\dots,u_{22}$. Table
    \ref{tab:linquot} shows all the linear quotients $J_{u_{i}}:u_i$, as
    described by \Cref{clm:G}.  After applying either the well-known formula
    \cite[Corollary 8.2.2]{MR2724673} or \Cref{SSSBetti}, we can obtain the
    following Betti diagram for $J$:
    \[
        \begin{array}{rrrrrrr}
            & 0 & 1 &  2 & 3 & 4 &5 \\
            \text{total}: &22 &75 &115 &94 &40 &7 \\
            1: &18 &56 & 79 &60 &24 &4 \\
            2: & 4 &19 & 36 &34 &16 &3
        \end{array}  
    \]
    Both the linear quotients and the Betti diagram agree with the outputs of \texttt{Macaulay2} \cite{M2}.
    \begin{table}[htbp]
        % TestSCM1(4,6,2)
        %i4 : BettiDia(4,6,2)
        %Here is the betti diagram of the dual of the polarization ideal, computed by Macaulay2
        %       0  1  2   3  4  5 6
        %total: 1 22 75 115 94 40 7
        %    0: 1  .  .   .  .  . .
        %    1: . 18 56  79 60 24 4
        %    2: .  4 19  36 34 16 3
        %
        %Here is the list of Lambda set of the dual ideal
        %{{1, 1}, {1, 2}, {2, 2, 2}}
        %
        %Here is a (maybe not complete) table of the computations by our formula
        %(i,d):beta_(i,i+d)
        %
        %(0, 0):0 (1, 0):0 (2, 0):0 (3, 0):0 (4, 0):0 (5, 0):0 (6, 0):0
        %(0, 1):0 (1, 1):0 (2, 1):0 (3, 1):0 (4, 1):0 (5, 1):0 (6, 1):0
        %(0, 2):18 (1, 2):56 (2, 2):79 (3, 2):60 (4, 2):24 (5, 2):4 (6, 2):0
        %(0, 3):4 (1, 3):19 (2, 3):36 (3, 3):34 (4, 3):16 (5, 3):3 (6, 3):0
        %(0, 4):0 (1, 4):0 (2, 4):0 (3, 4):0 (4, 4):0 (5, 4):0 (6, 4):0
        %(0, 5):0 (1, 5):0 (2, 5):0 (3, 5):0 (4, 5):0 (5, 5):0 (6, 5):0
        %(0, 6):0 (1, 6):0 (2, 6):0 (3, 6):0 (4, 6):0 (5, 6):0 (6, 6):0
        \caption{Linear quotients of a sparse symmetric shifted ideal}
        \label{tab:linquot}
        \begin{minipage}[t]{0.45\linewidth}
            \renewcommand{\arraystretch}{1.2}
            \centering
            {\tiny
                \begin{tabular}{cccc}
                    \hline
                    $i$ & $u_i$ & $J_{u_{i}}:u_i$ & $\max(u_i)$\\ \hline 
                    1 & $z_{1, 1} z_{2, 1}$ & - & $2$\\ 
                    2 & $z_{1, 1} z_{3, 1}$ & $\braket{z_{2, 1}}$ & $3$\\ 
                    3 & $z_{1, 1} z_{4, 1}$ & $\braket{z_{3, 1},z_{2, 1}}$ & $4$\\ 
                    4 & $z_{2, 1} z_{3, 1}$ & $\braket{z_{1, 1}}$ & $3$\\ 
                    5 & $z_{2, 1} z_{4, 1}$ & $\braket{z_{3, 1},z_{1, 1}}$ & $4$\\ 
                    6 & $z_{3, 1} z_{4, 1}$ & $\braket{z_{2, 1},z_{1, 1}}$ & $4$\\ 
                    7 & $z_{1, 1} z_{2, 2}$ & $\braket{z_{4, 1},z_{3, 1},z_{2, 1}}$ & $2$\\ 
                    8 & $z_{1, 1} z_{3, 2}$ & $\braket{z_{2, 2},z_{4, 1},z_{3, 1},z_{2, 1}}$ & $3$\\ 
                    9 & $z_{1, 1} z_{4, 2}$ & $\braket{z_{3, 2},z_{2, 2},z_{4, 1},z_{3, 1},z_{2, 1}}$ & $4$\\ 
                    10 & $z_{2, 1} z_{1, 2}$ & $\braket{z_{4, 1},z_{3, 1},z_{1, 1}}$ & $1$\\ 
                    11 & $z_{2, 1} z_{3, 2}$ & $\braket{z_{1, 2},z_{4, 1},z_{3, 1},z_{1, 1}}$ & $3$\\ \hline
                \end{tabular}}
        \end{minipage}
        \begin{minipage}[t]{0.45\textwidth}
            \renewcommand{\arraystretch}{1.2}
            \centering
            {\tiny
                \begin{tabular}{cccc}
                    \hline
                    $i$ & $u_i$ & $J_{u_{i}}:u_i$ & $\max(u_i)$\\ \hline
                    12 & $z_{2, 1} z_{4, 2}$ & $\braket{z_{3, 2},z_{1, 2},z_{4, 1},z_{3, 1},z_{1, 1}}$ & $4$\\ 
                    13 & $z_{3, 1} z_{1, 2}$ & $\braket{z_{4, 1},z_{2, 1},z_{1, 1}}$ & $1$\\ 
                    14 & $z_{3, 1} z_{2, 2}$ & $\braket{z_{1, 2},z_{4, 1},z_{2, 1},z_{1, 1}}$ & $2$\\ 
                    15 & $z_{3, 1} z_{4, 2}$ & $\braket{z_{2, 2},z_{1, 2},z_{4, 1},z_{2, 1},z_{1, 1}}$ & $4$\\ 
                    16 & $z_{4, 1} z_{1, 2}$ & $\braket{z_{3, 1},z_{2, 1},z_{1, 1}}$ & $1$\\ 
                    17 & $z_{4, 1} z_{2, 2}$ & $\braket{z_{1, 2},z_{3, 1},z_{2, 1},z_{1, 1}}$ & $2$\\ 
                    18 & $z_{4, 1} z_{3, 2}$ & $\braket{z_{2, 2},z_{1, 2},z_{3, 1},z_{2, 1},z_{1, 1}}$ & $3$\\ 
                    19 & $z_{1, 2} z_{2, 2} z_{3, 2}$ & $\braket{z_{4, 1},z_{3, 1},z_{2, 1},z_{1, 1}}$ & $3$\\ 
                    20 & $z_{1, 2} z_{2, 2} z_{4, 2}$ & $\braket{z_{3, 2},z_{4, 1},z_{3, 1},z_{2, 1},z_{1, 1}}$ & $4$\\ 
                    21 & $z_{1, 2} z_{3, 2} z_{4, 2}$ & $\braket{z_{2, 2},z_{4, 1},z_{3, 1},z_{2, 1},z_{1, 1}}$ & $4$\\ 
                    22 & $z_{2, 2} z_{3, 2} z_{4, 2}$ & $\braket{z_{1, 2},z_{4, 1},z_{3, 1},z_{2, 1},z_{1, 1}}$ & $4$\\ \hline
                \end{tabular}}
        \end{minipage}
    \end{table}  
\end{Example}

\subsection{General case}
Notice that in the standard monomial case, the big height of
$I_a(\calF_\bdz^b)$ is $h= bs-a+1$ by \cite[Remark 4.5]{arXiv:1912.04448}.
Thus, for the general case when $I=I_a(\calF^b)$, we need the following
genericness assumption.

\begin{Setting}
    \label{setting-general-case}
    With the assumptions in \Cref{setting-monomial-case}, we assume further
    that  $n\ge bs-a+1$ and $\calF=\{f_1,\dots,f_s\}$ is a set of
    $(bs-a+1)$-generic $d$-forms in $R=\KK[x_0,\dots,x_n]$.
\end{Setting}

From now on, we will fix a positive integer $m$ and investigate the symbolic
power $I_a(\calF^b)^{(m)}$.  It follows from \cite[Theoorem
4.6]{arXiv:1912.04448} that
\[
    I_a(\calF^b)^{(m)}
    = \bigcap_{c=c_0}^{h} I_{c,\calF}^{(m(\mu_a^0+b(c-c_0)))}.
\]
Inspired by the discussion in \cite[Subsection 4.2]{MR2184792}, we
simultaneously look at the ideals
\[
    K_j\coloneqq 
    \bigcap_{c=c_0}^{h-j}
    I_{c,\calF}^{(m(\mu_a^0+b(c-c_0)))}
    \qquad \text{and} \qquad
    K_j(\bdz)\coloneqq 
    \bigcap_{c=c_0}^{h-j}
    I_{c,\calF_{\bdz}}^{(m(\mu_a^0+b(c-c_0)))}
\]
for $j=0,1,\dots,h-c_0$. We may also write 
$K_{h-c_0+1}=R$ and 
$K_{h-c_0+1}(\bdz)=T$.
Thus, we have a filtration of submodules
\[
    0=\overline{I_a(\calF^b)^{(m)}}=\overline{K_0}\subseteq
    \overline{K_1}\subseteq \cdots \subseteq \overline{K_{h-c_0}}\subseteq
    R/I_a(\calF^b)^{(m)},
\]
where $\overline{\phantom m}$ means the corresponding image in
$R/I_a(\calF^b)^{(m)}$. And we have a similar filtration of submodules for
$T/I_a(\calF_{\bdz}^b)^{(m)}$.  Now,
$\overline{K_{j+1}(\bdz)}/\overline{K_j(\bdz)}\cong K_{j+1}(\bdz)/K_j(\bdz)$ is
either $0$ or Cohen--Macaulay of expected dimension for each $j$, by
\cite[Theorem A.4]{MR2184792} and \Cref{prop:mono-SCM}. It remains to check if
each $\overline{K_{j+1}}/\overline{K_j}$ still has the similar property.

\begin{Lemma}
    \label{lem:dim}
    \begin{enumerate}[a]
        \item If $\overline{K_{j+1}}\ne \overline{K_j}$, then
            $\dim(R)-\dim(\overline{K_{j+1}}/ \overline{K_j})=h-j$.
        \item Similarly, if $\overline{K_{j+1}(\bdz)}\ne \overline{K_j(\bdz)}$,
            then
            $\dim(T)-\dim(\overline{K_{j+1}(\bdz)}/\overline{K_j(\bdz)})=h-j$
            as well.
    \end{enumerate}
\end{Lemma}

\begin{proof}
    In the first case,
    \begin{align*}
        \overline{K_{j+1}}/\overline{K_j}&\cong
        {K_{j+1}}/{K_j}=K_{j+1}/(K_{j+1}\cap
        I_{c,\calF}^{(m(\mu_a^0+b(c-c_0)))})\\
        &\cong (K_{j+1}+ I_{c,\calF}^{(m(\mu_a^0+b(c-c_0)))})/
        I_{c,\calF}^{(m(\mu_a^0+b(c-c_0)))} \subset R/
        I_{c,\calF}^{(m(\mu_a^0+b(c-c_0)))}
    \end{align*}
    for $c=h-j$. Thus, $\Ass(\overline{K_{j+1}}/\overline{K_j})\subseteq
    \Ass(R/ I_{c,\calF}^{(m(\mu_a^0+b(c-c_0)))})$. But
    $R/I_{c,\calF}^{(m(\mu_a^0+b(c-c_0)))}$ is Cohen--Macaulay of dimension
    $\dim(R)-c$ by \cite[Corollary 3.7]{MR3683102}. This shows that
    $\dim(\overline{K_{j+1}}/\overline{K_j})=\dim(R)-c$ as well. The second
    case is similar.
\end{proof}

\begin{Lemma}
    \label{lem:resolution}
    Let $\varphi:T=\KK[z_1,\dots,z_s]\to R$ be the homomorphism induced by
    $z_i\mapsto f_i$ for each $i$.  Suppose that $\bdF_{\bullet}$ is a graded
    minimal free resolution of $T/K_j(\bdz)$ for
    $\calF_{\bdz}=\{z_1,\dots,z_s\}\subseteq T$. Then,
    $K_j=\varphi(K_j(\bdz))R$ and $\bdF_{\bullet}\otimes_{T} R$ is a graded
    minimal free resolution of $R/K_j$.  
\end{Lemma}

\begin{proof}
    The argument will be similar to the proof of \cite[Theorem 3.3]{MR3683102}.
    We prove by induction on $s-\hat{c}$, where $\hat{c}\coloneqq
    \max\{c:\text{$\calF$ is $c$-generic}\}$.  Obviously, $h\le \hat{c}\le
    s-1$. If $s-\hat{c}=1$, $\calF$ forms a regular sequence. Whence, $\varphi$
    is flat and the claim is clear by \cite[Remark 2.1]{arXiv:1912.04448} or
    \cite[Lemma 3.1]{MR3683102} and its proof.

    When $s-\hat{c}>1$, let $y$ be a new variable over $R$. For each $i$, let
    $f_i'$ be a general $d$-form in the ideal $(f_i,y)\subset R[y]$. Now,
    consider the new homomorphism $\gamma: T\to R[y]$ induced by $z_i\mapsto
    f_i'$. Notice that $\calF'\coloneqq\{f_1',\dots,f_s'\}$ is a set of
    $(\hat{c}+1)$-generic forms in $R[y]$. Whence, by induction,
    $\bdF_{\bullet}\otimes_T R[y]$ is a graded minimal free resolution of
    $R[y]/K_j'$ where $K_j'=\bigcap_{c=c_0}^{h-j}
    I_{c,\calF'}^{(m(\mu_a^0+b(c-c_0)))}$.  Meanwhile, we have the graded
    isomorphism
    \[
        R[y]/(K_j',y) \cong R/K_j.
    \]
    Thus, the last piece of the proof is to show that $y$ is a non-zero-divisor
    of $R[y]/K_j'$. 

    Suppose for contradiction that $y$ is not regular modulo the ideal $K_j'$,
    then $y \in \frakp$ for some associated prime $\frakp$. Notice that this is
    equivalent to saying that $y$ is not regular when modulo
    $I_{c,\calF'}^{(m(\mu_a^0+b(c-c_0)))}$ for some $c$ with $c_{0}\le c\le
    h-j$. But this is impossible by the proof of \cite[Theorem 3.3]{MR3683102}.
\end{proof}

Here is the main result of this section.

\begin{Theorem}
    %\label{thm:SCM}
    Assuming \Cref{setting-general-case}, the symbolic powers of
    $I_a(\calF^b)=I_a(f_1^b\cdots f_s^b)$ are all sequentially Cohen--Macaulay.
\end{Theorem}

\begin{proof}
    We will fix a positive integer $m$ and look at the symbolic power
    $I_a(\calF^b)^{(m)}$.  After applying the notations $K_j$ and $K_j(\bdz)$
    introduced before \Cref{lem:dim}, we have a short exact sequence
    \[
        0\to K_{j+1}(\bdz)/K_j(\bdz) \to T/K_j(\bdz) \to T/K_{j+1}(\bdz) \to 0.
    \]
    We may assume that $K_j\ne K_{j+1}$. Hence, by \Cref{lem:resolution},
    $K_j(\bdz)\ne K_{j+1}(\bdz)$. Furthermore,
    \[
        \Tor_k^T(T/K_j(\bdz),R)=0=\Tor_k^T(T/K_{j+1}(\bdz),R) 
    \]
    for $k\ge 1$. If we take the well-known induced long exact sequence, then
    we obtain $\Tor_k^T(K_{j+1}(\bdz)/K_j(\bdz),R)=0$ as well for all $k\ge 1$.
    Furthermore, since $(T/K_j(\bdz))\otimes R\cong R/K_j$,
    $(T/K_{j+1}(\bdz))\otimes R\cong R/K_{j+1}$ and
    $\Tor_1(T/K_{j+1}(\bdz),R)=0$, we see immediately that
    $(K_{j+1}(\bdz)/K_j(\bdz))\otimes R\cong K_{j+1}/K_j$.  Therefore, if
    $\bdG_\bullet$ is a graded minimal free resolution of
    $K_{j+1}(\bdz)/K_j(\bdz)$, then $\bdG_\bullet\otimes R$ will be a graded
    minimal free resolution of $K_{j+1}/K_j$. In particular,
    $\pd(K_{j+1}/K_j)=\pd(K_{j+1}(\bdz)/K_j(\bdz))$.

    Since $K_{j+1}(\bdz)/K_j(\bdz)$ is Cohen--Macaulay, it follows from
    \Cref{lem:dim} and the Auslander--Buchsbaum formula that $K_{j+1}/K_{j}$ is
    also Cohen--Macaulay of the expected dimension. Therefore,
    $I_a(\calF^b)^{(m)}$ is sequentially Cohen--Macaulay.
\end{proof}

\section{Harbourne--Huneke containment problem for the uniform $a$-fold product ideal}
Recall that for any homogeneous ideal $I$ in a standard graded ring,
$\alpha(I)$ stands for the least degree of a homogeneous polynomial in this
ideal. This final section is mainly concerned with the following two questions
involving this invariant, which were raised for \emph{radical} homogeneous
ideals in \cite[Section 3]{arXiv:2009.05022}.

\begin{Question}\label{Delike}
    [Demailly-like bound]
    Let $R$ be a polynomial ring over a field $\KK$ and let $I$ be a
    homogeneous ideal of big height $h$ in $R$. Does the inequality
    \[
        \frac{\alpha(I^{(n)})}{n}\ge \frac{\alpha(I^{(m)})+h-1}{m+h-1}
    \]
    hold for all $n, m \ge 1$? 
\end{Question}

Notice that the \emph{Waldschmidt constant}
\[
    \widehat{\alpha}(I)\coloneqq \lim_n \frac{\alpha(I^{(n)})}{n}
\]
actually agrees with $\inf_n \frac{\alpha(I^{n})}{n}$. Thus, the above inequality
is simply
\[
    \widehat{\alpha}(I)\ge \frac{\alpha(I^{(m)})+h-1}{m+h-1}.
\]

\begin{Question}\label{HHContain}
    [General version of the Harbourne--Huneke containment for radical homogeneous ideal]
    Let $R$ be a polynomial ring over a field $\KK$ with the maximal
    homogeneous ideal $\frakm$. If $I$ is a homogeneous ideal of big height $h$
    in $R$, does the inequality
    \[
        I^{(\ell(h+m-1))}\subseteq \frakm^{\ell(h-1)}(I^{(m)})^\ell
    \]
    holds for all $\ell, m\ge 1$?
\end{Question}

We will drop the radical requirement and study these questions for
$I=I_a(\calF^b)$ under the assumptions in \Cref{setting-general-case}.  It is
not difficult to see that a positive answer to \Cref{HHContain} would imply an
affirmative answer to \Cref{Delike}. Thus, we will focus on \Cref{HHContain}.
Actually, like \cite[Theorem 3.6]{arXiv:2009.05022}, a stronger result can be
achieved. Notice that the big height of the ideal $I_a(\calF^b)$ is $h=bs-a+1$
while the height is the positive integer $c_0$ by \cite[Remark
4.5]{arXiv:1912.04448}.  Furthermore, we have seen in the proof of
\cite[Proposition 4.7]{arXiv:1912.04448} that $c_0=h$ precisely when $b=1$ or
$a=bs$.  The $b=1$ case is just the standard star configuration case studied in
\cite{arXiv:2009.05022}. Meanwhile, the $a=bs$ case is the trivial one when the
ideal is actually principal.  Thus, in the following, we may assume that
$c_0<h$. Consequently, $b>1$, $a<bs$ and $h-1>0$.

We present two technical lemmas first.

\begin{Lemma}
    \label{lem-eqn-8}
    We have
    %\begin{equation}
    $(h-1)(\mu_a^0+b(c-c_0))\ge c$
    %\label{eqn-8}
    %\end{equation}
    for each $c$ with $c_0\le c\le h$.
\end{Lemma}

\begin{proof}
    Notice that the coefficients of $c$ are $(h-1)b$  and $1$ on the two sides. Since $(h-1)b\ge 1$, it suffices to show the inequality for $c=c_0$, namely, to show 
    \begin{equation}
        (h-1)\mu_a^0 \ge c_0.
        \label{eqn-4}
    \end{equation} 
    After putting back the definitions of $h$,  $\mu_a^0$ and $c_0$, this inequality is actually
    \[
        (bs-a)\left(a-b\floor{\frac{a-1}{b}}\right)\ge s-\floor{\frac{a-1}{b}}.
    \]
    And it can be simplified into
    \[
        (bs-a)a -s \ge \left( (bs-a)b-1 \right)\floor{\frac{a-1}{b}}.
    \]
    Since $bs>a$ and $b> 1$, we have $(bs-a)b-1>0$. Thus, it suffices to show that
    \[
        (bs-a)a -s \ge \left( (bs-a)b-1 \right)\cdot \frac{a-1}{b},
    \]
    which is equivalent to saying that
    \[
        (b-1)(bs-a)\ge 1.
    \]
    However, the last inequality is clear. Thus, we have successfully shown the inequality \eqref{eqn-4}. And this completes our proof for \Cref{lem-eqn-8}.
\end{proof}

In the following, for each positive integer $c$ such that $c_0\le c\le h$, we write $\widehat{\mu}_c$ for the positive integer $\mu_a^0+b(c-c_0)$.  

\begin{Lemma}
    \label{lem-eqn-6}
    Let $k,\ell$, and $m$ be positive integers.
    Then, we have 
    %\begin{equation}
    \[
        (\ell(h+m-1)-h+k)\widehat{\mu}_c-c(\ell-1)\ge \ell m \widehat{\mu}_c
    \]
    %\label{eqn-6}
    %\end{equation}
    for each $c$ with $c_0\le c\le h$.
\end{Lemma}

\begin{proof}
    This artificially irregular inequality can be quickly modified into
    %\[
    %    \left( \ell(h-1)-h+k \right)\widehat{\mu}_c\ge c(\ell-1),
    %\]
    %namely,
    \begin{equation}
        \ell\left( (h-1)\widehat{\mu}_c-c \right)\ge (h-k)\widehat{\mu}_c-c.
        \label{eqn-7}
    \end{equation}
    Since $(h-1)\widehat{\mu}_c\ge c$ by \Cref{lem-eqn-8}, we only need to show \eqref{eqn-7} for $\ell=1$, i.e., to show that
    \[
        (h-1)\widehat{\mu}_c-c \ge (h-k)\widehat{\mu}_c-c.
    \]
    However, this holds automatically, since both $k$ and $\widehat{\mu}_c$ are positive integers.
\end{proof}

%In the following,  we will write $m'\coloneqq \ell(h+m-1)-h+k$. 

Now, we are ready for the main theorem of this section.
\begin{Theorem}
    %\label{thmHH}
    Under the assumptions in \Cref{setting-general-case}, the ideal 
    $I=I_a(\calF^b)$ satisfies
    \begin{equation}
        I^{(\ell(h+m-1)-h+k)}\subseteq \frakm^{d((\ell-1)(h-1)+k-1)\widehat{\mu}_h}(I^{(m)})^\ell
        \label{eqn-5}
    \end{equation}
    for all positive integers $k,\ell$ and $m$.
\end{Theorem}

Note that $\widehat{\mu}_h$ above is simply $\mu_a^0+b(h-c_0)=bh-(bs-a)$.

\begin{proof}
    Firstly, we consider the standard monomial case when
    $I=I_a(\calF_{\bdz}^b)$ for $\calF_\bdz=\{z_1,\dots,z_s\}$ in
    $T=\KK[z_1,\dots,z_s]$.  Within this ring, by abuse of notation, we will
    write $\bdz^\bdmu \coloneqq z_1^{\mu_1}\cdots z_s^{\mu_s}$ for any tuple
    $\bdmu=(\mu_1,\dots,\mu_s)\in \ZZ_{\ge0}^s$. As usual, $|\bdmu|$ is simply
    $\mu_1+\cdots+\mu_s$, and $\bdmu_{\le c}$ is short for the sub-tuple
    $(\mu_1,\dots,\mu_c)$.  Since the permutation group $\frakS_s$ has a
    natural action on $T$, for any monomial ideal $K$ that is invariant under this
    action, we will introduce
    \[
        \calP(K)\coloneqq\Set{\bdlambda\text{ is a partition of length $s$}: \bdz^{\bdlambda}\in K} 
    \]
    and, by abuse of notation, 
    \[
        \Lambda(K)\coloneqq\Set{\bdlambda:\bdz^{\bdlambda} \text{ is a minimal monomial generator of $K$}}.
    \]

    Take an arbitrary partition $\bdlambda=(\lambda_1,\dots,\lambda_s)\in
    \calP(I_a(\calF_{\bdz}^b)^{(m')})$ where $m'\coloneqq \ell(h+m-1)-h+k$.
    Due to the decomposition in  \eqref{eqn:IaDecomp},
    %\cite[Theorem 4.6]{arXiv:1912.04448}, 
    this amounts to saying that
    \[
        \bdlambda\in \calP(I_{c,\calF_{\bdz}}^{(m'\widehat{\mu}_c)})\qquad
        \text{for each $c$ with $c_0\le c\le h$,}
    \]
    namely,
    \[
        |\bdlambda_{\le c}|\ge m'\widehat{\mu}_c \qquad
        \text{for each $c$ with $c_0\le c\le h$}
    \]
    by \cite[Proposition 4.1]{MR4108337}.

    In the following, for each $j\le h$, we write $\lambda_j'\coloneqq
    \floor{\frac{\lambda_j}{\ell}}$ and $\widehat{\lambda}_j\coloneqq
    \lambda_j-\ell\lambda_j'$. Since $\widehat{\lambda}_j\le \ell-1$ for such
    $j$, we have $\sum_{j\le c}\widehat{\lambda}_j\le c(\ell-1)$. Thus,
    \[
        \sum_{j\le c}\ell\lambda_j' =\sum_{j\le c}(\lambda_j-\widehat{\lambda}_j)=\sum_{j\le c}\lambda_j -\sum_{j\le c}\widehat{\lambda}_j\\
        \ge m'\widehat{\mu}_c-c(\ell-1).
    \]
    As a quick consequence of \Cref{lem-eqn-6}, we have
    \[
        \sum_{j\le c}\lambda_j'\ge m\widehat{\mu}_c
    \]
    for each such $c$.  It follows again from \cite[Proposition 4.1]{MR4108337} and the decomposition in  \eqref{eqn:IaDecomp}
    %\cite[Theorem 4.6]{arXiv:1912.04448}  
    that
    \[
        \bdlambda'\coloneqq (\lambda_1',\lambda_2',\dots,\lambda_h',\lambda_h',\dots,\lambda_h')\in \calP(I_a(\calF_{\bdz}^b)^{(m)}).
    \]
    Therefore, we can find suitable partition
    $\bdlambda''=(\lambda_1'',\dots,\lambda_s'')\in
    \Lambda(I_a(\calF_{\bdz}^b)^{(m)})$ such that $\bdlambda''\le \bdlambda'$
    componentwise.  In the part (b) of the proof of \cite[Proposition
    4.7]{arXiv:1912.04448}, we have shown that $|\bdlambda_{\le
        h}''|=m\widehat{\mu}_{h}$.  Furthermore, $\lambda_j''=\lambda_h''$
    whenever $h\le j\le s$. Now,  $\ell\bdlambda''\le \bdlambda$ componentwise.
    And it remains to estimate $|\bdlambda|-\ell|\bdlambda''|$. It is clear
    that
    \begin{align*}
        |\bdlambda|-\ell|\bdlambda''|&=\sum_{j=1}^s (\bdlambda_j-\ell\bdlambda_j'')
        \ge \sum_{j=1}^h (\bdlambda_j-\ell\bdlambda_j'')
        =|\bdlambda_{\le h}|-\ell |\bdlambda_{\le h}''|\\
        &\ge m'\widehat{\mu}_h-\ell m \widehat{\mu}_h 
        =(\ell(h+m-1)-h+k-\ell m)\widehat{\mu}_h\\
        &= \left( (\ell-1)(h-1)+k-1 \right)\widehat{\mu}_h.
    \end{align*}

    Consequently, for $I=I_a(\calF_{\bdz}^b)$, we have
    \begin{equation}
        I^{(\ell(h+m-1)-h+k)}\subseteq \frakm^{\left((\ell-1)(h-1)+k-1 \right)\widehat{\mu}_h}(I^{(m)})^\ell,
        \label{eqn-HH-2}
    \end{equation}
    %Notice that in the standard monomial star configuration case, $b=1$. Whence, $h=c_0$ and $\widehat{\mu}_h=1$.
    %As $\widehat{\mu}_h\ge 1$, this shows 
    establishing the containment \eqref{eqn-5} in this standard monomial case.

    In the final step, we consider the general uniform star configuration case,
    namely when $I=I_a(\calF^b)$. Notice that both $I^{(\ell(h+m-1)-h+k)}$ and
    $I^{(m)}$ are generated by monomials in the forms in $\calF$;
    cf.~\cite[Proposition 4.11]{arXiv:1912.04448}.  Thus, it follows directly
    from \eqref{eqn-HH-2} that
    \[
        I^{(\ell(h+m-1)-h+k)}\subseteq \frakm^{d\left( (\ell-1)(h-1)+k-1 \right)\widehat{\mu}_h}(I^{(m)})^\ell,
    \]
    where $d$ is the common degree of the forms in $\calF$. 
\end{proof}

The following result generalizes \cite[Corollary 3.7]{arXiv:2009.05022}.

\begin{Corollary}
    The Harbourne–Huneke type containment \textup{(}\Cref{HHContain}\textup{)} and the Demailly-like bound \textup{(}\Cref{Delike}\textup{)} hold for the ideal $I=I_a(\calF^b)$ under the assumptions in \Cref{setting-general-case}.   
\end{Corollary}

\begin{proof}
    The Harbourne–Huneke type containment follows from \eqref{eqn-5}, since both $d$ and $\widehat{\mu}_h$ are positive integers, while $k$ can be specially taken to be $h$. By applying the identical arguments in the proof of \cite[Corollary 3.7]{arXiv:2009.05022}, we see that the Demailly-like bound also holds.
\end{proof}

\begin{acknowledgment*}
    The authors are grateful to the software system \texttt{Macaulay2} \cite{M2}, for serving as an excellent source of inspiration.
    The second author is partially supported by the ``Anhui Initiative in Quantum Information Technologies'' (No.~AHY150200).
\end{acknowledgment*}

\end{document}